\documentstyle[12pt]{article}
\begin{document}
\newpage
\pagestyle{empty}
\setcounter{page}{0}
%
\newfont{\twelvemsb}{msbm10 scaled\magstep1}
\newfont{\eightmsb}{msbm8}
\newfont{\sixmsb}{msbm6}
\newfam\msbfam
\textfont\msbfam=\twelvemsb
\scriptfont\msbfam=\eightmsb
\scriptscriptfont\msbfam=\sixmsb
\catcode`\@=11
\def\Bbb{\ifmmode\let\next\Bbb@\else
  \def\next{\errmessage{Use \string\Bbb\space only in math mode}}\fi\next}
\def\Bbb@#1{{\Bbb@@{#1}}}
\def\Bbb@@#1{\fam\msbfam#1}
\newfont{\twelvegoth}{eufm10 scaled\magstep1}
\newfont{\tengoth}{eufm10}
\newfont{\eightgoth}{eufm8}
\newfont{\sixgoth}{eufm6}
\newfam\gothfam
\textfont\gothfam=\twelvegoth
\scriptfont\gothfam=\eightgoth
\scriptscriptfont\gothfam=\sixgoth
\def\frak{\frak@}
\def\frak@#1{{\fam\gothfam{{#1}}}}
\def\frak@@#1{\fam\gothfam#1}
\catcode`@=12

%
%
%
\def\CC{{\Bbb C}}
\def\NN{{\Bbb N}}
\def\QQ{{\Bbb Q}}
\def\RR{{\Bbb R}}
\def\ZZ{{\Bbb Z}}
\def\cA{{\cal A}}          \def\cB{{\cal B}}          \def\cC{{\cal C}}
\def\cD{{\cal D}}          \def\cE{{\cal E}}          \def\cF{{\cal F}}
\def\cG{{\cal G}}          \def\cH{{\cal H}}          \def\cI{{\cal I}}
\def\cJ{{\cal J}}          \def\cK{{\cal K}}          \def\cL{{\cal L}} 
\def\cM{{\cal M}}          \def\cN{{\cal N}}          \def\cO{{\cal O}}
\def\cP{{\cal P}}          \def\cQ{{\cal Q}}          \def\cR{{\cal R}} 
\def\cS{{\cal S}}          \def\cT{{\cal T}}          \def\cU{{\cal U}}
\def\cV{{\cal V}}          \def\cW{{\cal W}}          \def\cX{{\cal X}}
\def\cY{{\cal Y}}          \def\cZ{{\cal Z}}
\def\qed{\hfill \rule{5pt}{5pt}}
\newtheorem{lemma}{Lemma}
\newtheorem{prop}{Proposition}
\newtheorem{theo}{Theorem}
\newenvironment{result}{\vspace{.2cm} \em}{\vspace{.2cm}}

$$
\;
$$
\rightline{CPTH-S538.0697}
\rightline{q-alg/9706033}
\rightline{June 96}

\vfill
\vfill
\begin{center}

  {\LARGE {\bf {\sf General Construction of Nonstandard $R_h$-matrices 
as Contraction Limits of $R_{q}$-matrices}}}
 \\[2cm]

\smallskip 

{\large B. Abdesselam\footnote{abdess@orphee.polytechnique.fr},
A. Chakrabarti\footnote{chakra@orphee.polytechnique.fr} and 
R. Chakrabarti\footnote{Permanent address: Department of Theoretical Physics, University of Madras, Guindy Campus, Madras-600025, India}}

\smallskip 

\smallskip 

\smallskip 

{\em  \footnote{Laboratoire Propre du CNRS UPR A.0014}Centre de Physique 
Th\'eorique, Ecole Polytechnique, \\
91128 Palaiseau Cedex, France.}

\end{center}

\vfill

\begin{abstract}
\noindent A class of transformations of $R_q$-matrices is introduced such that
the $q\rightarrow 1$ limit gives explicit nonstandard $R_{h}$-matrices.
The transformation matrix is singular itself at $q\rightarrow 1$ limit. 
For the transformed matrix, the singularities, however, cancel yielding a 
well-defined construction. Our method 
can be implemented systematically for $R$-matrices of all dimensions 
and not only for $sl(2)$ but also for algebras of higher dimensions. Explicit 
constructions are presented starting with ${\cal U}_q(sl(2))$ and ${\cal U}_q
(sl(3))$, while choosing $R_q$ for {\em (fund. rep.)$\otimes$(arbitrary 
irrep.)}. The treatment for the general case and various perspectives 
are indicated. Our method yields nonstandard deformations along with a 
nonlinear map of the $h$-Borel subalgebra on the corresponding classical 
Borel subalgebra. For ${\cal U}_h(sl(2))$ this map is extended to the whole 
algebra and compared with another one proposed by us previously.      
\end{abstract}

\vfill
\vfill

\newpage
\pagestyle{plain}

The $R$-matrices for the fundamental representations of the nonstandard 
$h$-deformations of $sl(2)$ and $so(4)(\simeq sl(2)\otimes sl(2))$ were 
obtained [1,2] through a specific contraction of the corresponding 
$q$-deformed $R$-matrices. A similarity transformation of the $4\times 4$ 
$R_q$-matrix for the fundamental representation of ${\cal U}_q(sl(2))$ was 
performed using a transforming matrix singular itself at the $q\rightarrow 1$ 
limit, but in such way that all singularities cancel out for the
transformed $R$-matrix giving the finite nonstandard $R_h$-matrix.
Following the previous practice [1,2], we refer to this combined process 
of similarity transformation and subsequent cancellation of singularities at 
the $q\rightarrow 1$ limit as contraction procedure. This 
technique was generalized to higher dimensional algebras [3] considering again 
the $N^2\times N^2$ dimensional $R$-matrices for the fundamental 
representations of 
$q$-deformed $sl(N)$, for example. Other relevant references can be found 
in [1-3].    

\smallskip

We present here an operatorial generalization of this approach directly 
applicable to $R$-matrices of arbitrary dimensions. For brevity and 
simplicity we start with $({1\over 2}\otimes j)$ representation i.e. 
$2(2j+1) \times 2(2j+1)$ dimensional $R_q$-matrix for ${\cal U}_q(sl(2))$. 
Then we 
will indicate possible generalizations in different directions, using 
${\cal U}_q(sl(3))$ as a particular example. The universal ${\cal R}$-matrix 
for ${\cal U}_h(sl(2))$ has been given a particularly convenient form [4,5]. 
For $({1\over 2}\otimes j)$ representation this reduces to 
\begin{eqnarray}
R_h=\pmatrix{e^{hX}&& -hH+{h\over 2}(e^{hX}-e^{-hX}) \cr
&& \cr
0 && e^{-hX}\cr}.
\end{eqnarray}
Here $(H,X)$ are the generators of the Borel subalgebra of ${\cal U}_h(sl(2))$ 
satisfying 
\begin{eqnarray}
[H,X]= 2 {\sinh hX \over h}.
\end{eqnarray}
Absent from the upper triangular form (1) and indeed form the universal
${\cal R}_h$-matrix is the generator $Y$ completing the ${\cal U}_h(sl(2))$
algebra, namely satisfying 
\begin{eqnarray}
  && [H,Y]=-Y(\cosh hX)-(\cosh hX)Y, \qquad\qquad  
   [X,Y]=H. 
\end{eqnarray} 
We will show how (1) can be obtained, directly and for arbitrary $j$, from 
the corresponding $R_q$-matrix for (${1\over 2}\otimes j$) representation
given by ( see, for example [6])
\begin{eqnarray}
R_q=\pmatrix{q^{{\cal  J}_0/2}&& q^{1/2} (1-q^{-2}) {\cal J}_- \cr
&& \cr
0 && q^{-{\cal J}_0/2}\cr}.
\end{eqnarray}
Here the generators of ${\cal U}_q(sl(2))$ are denoted by 
$(q^{\pm {\cal J}_0},
{\cal J}_{\pm})$ satisfying the standard relations 
\begin{eqnarray}
  && q^{{\cal J}_{0}}{\cal J}_{\pm} = {\cal J}_{\pm}q^{{\cal J}_{0}\pm 2},
\qquad\qquad
 [{\cal J}_{+},{\cal J}_{-}] = {q^{{\cal J}_{0}}-q^{-{\cal J}_{0}}
\over q-q^{-1}}\equiv [{\cal J}_{0}].  
\end{eqnarray}
Throughout we consider generic $q$, excluding roots of unity.

\smallskip

For the purpose of transforming the $R_q$ matrix in (4), we now consider 
a $q$-deformed exponential operator: 
 \begin{eqnarray}
E_q(\eta {\cal J}_{+})=\sum_{n= 0}^{\infty} {(\eta{\cal J}_{+})^n\over [n]!},
\end{eqnarray}
where 
 \begin{eqnarray}
[n]\equiv {q^{n}-q^{-n} \over q-q^{-1}}, \qquad\qquad [n]!=[n][n-1]\cdots
[1],\qquad \qquad [0]!\equiv 1.\nonumber 
\end{eqnarray}
We choose, for an arbitrary finite constant $h$, the parameter $\eta$ as
\begin{eqnarray}
\eta={h\over q-1}. 
\end{eqnarray}
We emphasize that though the deformed exponential $E_q(x)$ defined in (6) has 
no convenient simple expression for its inverse (comparable to the standard 
$q$-exponential [7] satisfying $(exp_{q}(x))^{-1}= exp_{q^{-1}}(-x)$), this 
is precisely what is needed for obtaining non-singular limiting forms for the 
$R$-matrix elements for arbitrary representations and other interesting 
properties. For any given value of $j$, the series (6) may be terminated 
after setting ${\cal J}_{+}^{2j+1}=0$; but, we proceed quite generally as 
follows. Defining  
\begin{eqnarray}
&& {\cal T}_{(\alpha)}=(E_q(\eta {\cal J}_{+}))^{-1}E_q(q^{\alpha}\eta 
{\cal J}_{+}) 
\end{eqnarray}
with ${\cal T}_{(0)}=1$, we obtain  
\begin{eqnarray}
&& (E_q(\eta {\cal J}_{+}))^{-1}q^{\alpha {\cal J}_0\over 2}E_q(\eta 
{\cal J}_{+})={\cal T}_{(\alpha)}q^{\alpha {\cal J}_0\over 2} . 
\end{eqnarray} 
For transforming the $R_q$ matrix in (4), the operators ${\cal T}_{(\pm 1)}$ 
are of particular importance. We will be concerned with simple rational 
values of $\alpha$. For later use, we note the identity 
\begin{eqnarray}
&& (E_q(\eta {\cal J}_{+}))^{-1}q^{{(\alpha+\beta)\over 2}{\cal J}_0 }E_q(\eta 
{\cal J}_{+})=\biggl((E_q(\eta{\cal J}_{+}))^{-1}q^{{\alpha \over 2}{\cal J}_0}
E_q(\eta {\cal J}_{+})\biggr)\biggl((E_q(\eta {\cal J}_{+}))^{-1}
q^{{\beta \over 2} {\cal J}_0}E_q(\eta {\cal J}_{+})\biggr),\nonumber 
\end{eqnarray} 
which in notation (8) reads
\begin{eqnarray}
{\cal T}_{(\alpha+\beta)}q^{{\alpha +\beta\over 2}
{\cal J}_0}=({\cal T}_{(\alpha)}
q^{{\alpha \over 2}{\cal J}_0 })({\cal T}_{(\beta)}
q^{{\beta \over 2}{\cal J}_0 }) .
 \end{eqnarray} 
Moreover, using the identity
 \begin{eqnarray}
&& [{\cal J}_{+}^{n}, {\cal J}_{-}]={[n]\over q-q^{-1}}\biggl(
q^{{\cal J}_0/2}{\cal J}_{+}^{n-1}q^{{\cal J}_0/2}-q^{-{\cal J}_0/2}{\cal J}_{+}^{n-1}q^{-{\cal J}_0/2}\biggr)
\end{eqnarray}
the following commutator is obtained
 \begin{eqnarray}
&&  [E_q(\eta{\cal J}_{+}), {\cal J}_{-}]={\eta\over q-q^{-1}}\biggl(
q^{{\cal J}_0/2}E_q(\eta{\cal J}_{+})q^{{\cal J}_0/2}-q^{-{\cal J}_0/2}E_q(\eta{\cal J}_{+}) q^{-{\cal J}_0/2}
\biggr)\nonumber \\
&&\phantom{[E_q(\eta{\cal J}_{+}), {\cal J}_{-}]}={\eta\over q-q^{-1}}\biggl(
E_q(q\eta{\cal J}_{+})q^{{\cal J}_0}-E_q(q^{-1}\eta{\cal J}_{+})q^{-{\cal J}_0}
\biggr),
\end{eqnarray}
which, in turn, leads to
\begin{eqnarray}
&&(E_q(\eta {\cal J}_{+}))^{-1}{\cal J}_{-} E_q(\eta 
{\cal J}_{+})= -{\eta \over q-q^{-1}} \biggl({\cal T}_{(1)}q^{{\cal J}_0}-
{\cal T}_{(-1)}q^{-{\cal J}_0}\biggr)+{\cal J}_{-}.
\end{eqnarray}
Evaluating term by term, the $q\rightarrow 1$ limits of ${\cal T}_{(\pm 1)}$ 
are found to be {\em finite} and of the form 
\begin{eqnarray}
\lim_{q\rightarrow 1} {\cal T}_{(\pm 1)}=T_{(\pm 1)}=\sum_{n=0}^{\infty}
c_{n}^{(\pm)}(hJ_{+})^n,
\end{eqnarray}
where ($J_{0}, J_{\pm}$) are the generators of the classical $sl(2)$ algebra:
\begin{eqnarray}
[J_{0}, J_{\pm}]=\pm J_{\pm}, \qquad\qquad [J_{+},J_{-}]=J_{0}.
 \end{eqnarray}
The first few coefficients $\{c_{n}^{(\pm)}\;|\;n=0,1,2,...\}$ in (14) read  
\begin{eqnarray}
&& c_{0}^{(\pm)}=1,\qquad c_{1}^{(\pm)}=\pm 1,\qquad 
c_{2}^{(\pm)}={1\over 2},\qquad
 c_{3}^{(\pm)}=0,\qquad c_{4}^{(\pm)}=- {1\over 8},\qquad 
c_{5}^{(\pm)}=0 . \qquad
\end{eqnarray}
The series $(E_q(x))^{-1}$ can be constructed systematically upto any given 
order in $\eta$; and, consequently, the operators $T_{(\pm 1)}$ are obtained. 
Along these lines a program in MAPLE for evaluation of the coefficients 
$c_n^{(\pm)}$ is easy to make. If the limits are indeed finite, then from 
(10) it is evident that 
  \begin{eqnarray}
T_{(\alpha)}=(T_{(1)})^{\alpha},
\end{eqnarray} 
where $T_{(\alpha)}=\lim_{q\rightarrow 1}{\cal T}_{(\alpha)}$. 
Henceforth we write 
 \begin{eqnarray}
T_{(1)}=T.
\end{eqnarray}
The result (16) 
suggests the following derivation of the {\em closed} form of $T$. To this 
end, we left and right multiply the second commutation relation in (5) by 
$(E_q(\eta{\cal J}_0))^{-1}$ and $E_q(\eta{\cal J}_0)$ respectively:
\begin{eqnarray}
(E_q(\eta {\cal J}_{+}))^{-1}(q^{{\cal J}_0}-q^{-{\cal J}_0})
E_q(\eta {\cal J}_{+})=(q-q^{-1})
(E_q(\eta {\cal J}_{+}))^{-1}
[{\cal J}_{+},{\cal J}_{-}]
E_q(\eta {\cal J}_{+}).
\end{eqnarray} 
Using (7), (9) and (13), we obtain
 \begin{eqnarray}
&& {\cal T}_{(2)}q^{{\cal J}_0}-{\cal T}_{(-2)}q^{-{\cal J}_0}=-{h \over q-1}
\biggl({\cal T}_{(1)}
({\cal J}_{+}q^{{\cal J}_0}-q^{{\cal J}_0}{\cal J}_{+})-
{\cal T}_{(-1)}({\cal J}_{+}q^{-{\cal J}_0}-q^{-{\cal J}_0}{\cal J}_{+})\biggr)\nonumber \\
&&\phantom{{\cal T}_{(2)}q^{{\cal J}_0}-{\cal T}_{(-2)}q^{-{\cal J}_0}= } +  (q-q^{-1})
[{\cal J}_{+},{\cal J}_{-}]\nonumber \\
&& \phantom{{\cal T}_{(2)}q^{{\cal J}_0}-
{\cal T}_{(-2)}q^{-{\cal J}_0}}=h(q+1) 
\biggl({\cal T}_{(1)}{\cal J}_{+}q^{{\cal J}_0}+
q^{-2}{\cal T}_{(-1)}{\cal J}_{+}q^{-{\cal J}_0}\biggr)
+(q^{{\cal J}_0}-q^{-{\cal J}_0}).
\end{eqnarray}
Using (17), we now obtain the following equation for $T$ at the limit 
$q\rightarrow 1$:
\begin{eqnarray}
T^2-T^{-2}=(T+T^{-1})(2hJ_{+}), \nonumber 
\end{eqnarray}
which after a factorization yields
\begin{eqnarray}
T-T^{-1}=2hJ_{+}.
\end{eqnarray}
The quadratic relation (21) in $T$ is now solved:
\begin{eqnarray}
T^{\pm 1}=\pm hJ_{+}+\biggl(1+(hJ_{+})^2\biggr)^{1/2}.
\end{eqnarray}
This is our crucial result. 

\smallskip

With all these results now in hand we go back to $R_{q}$ in (4). We choose 
the transformation matrix as $G=g_{1/2}\otimes g$, where 
$g=E_q(\eta {\cal J}_{+})$ and $g_{1/2}\equiv g|_{j={1\over 2}}$. The 
similarity transformation now yields  
 \begin{eqnarray}
G^{-1}R_qG=\pmatrix{g^{-1}q^{{\cal J}_0/2}g&& 
\eta g^{-1}(q^{{\cal J}_0/2}-q^{-{\cal J}_0/2} )g +
q^{-1/2}(q-q^{-1})g^{-1}{\cal J}_-g \cr
&& \cr
0 && g^{-1}q^{-{\cal J}_0/2}g\cr}.
\end{eqnarray}
Implementing the relevant results from (9) to (23), we now obtain
\begin{eqnarray}
&& R_h=\lim_{q\rightarrow 1}(G^{-1}R_{q}G)=
\pmatrix{T&& -{h\over 2}(T+T^{-1})J_0+{h\over 2}(T-T^{-1}) \cr
&& \cr
0 && T^{-1}\cr} \\
&& \phantom{{\cal R}_h=\lim_{q\rightarrow 1}(G^{-1}R_{q}G)}=
\pmatrix{e^{h{\hat X}}&& -h{\hat H}+{h\over 2}(e^{h{\hat X}}
-e^{-h{\hat X}}) \cr
&& \cr
0 && e^{-h{\hat X}}\cr},
\end{eqnarray}
where we have defined 
\begin{eqnarray}
&& e^{\pm h{\hat X}}=T^{\pm 1}=\pm hJ_{+}+(1+(hJ_{+})^2)^{1/2}  \\
&& {\hat H}={1\over 2}(T+T^{-1})J_0=(1+(hJ_{+})^2)^{1/2}J_0.
\end{eqnarray}
It is easily to verify that 
\begin{eqnarray}
  && [{\hat H},T^{\pm 1}]=T^{\pm 2}-1\qquad \Rightarrow\qquad
  [{\hat H},{\hat X}]={2\over h}\sinh(h{\hat X}). 
\end{eqnarray} 
Comparing (1), (2) with (25), (28) respectively we see that the contraction 
sheme, which comprised of our 
transformation and limiting procedure has furnished the $R_{h}$-matrix
along with a {\em nonlinear map of the Borel subalgebra of ${\cal U}_h(sl(2))$
 generated by $(H,X)$ on the classical one generated by $(J_0,J_{+})$}.
Indeed, defining 
\begin{eqnarray}
{\hat Y}=J_{-}-{h^2\over 4}J_{+}(J_0^2-1)
\end{eqnarray} 
we complete the ${\cal U}_h(sl(2))$ algebra and show that
\begin{eqnarray}
&& [{\hat H},{\hat Y}]=-({\hat Y}\cosh(h{\hat X})
+\cosh(h{\hat X}){\hat Y}), \qquad\qquad
 [{\hat X},{\hat Y}]={\hat H}.
\end{eqnarray} 
The triplet $({\hat H},{\hat X},{\hat Y})$ may be looked as a particular 
realization of the ${\cal U}_h(sl(2))$ generators $(H,X,Y)$. 
We have, therefore developed an invertible nonlinear map of $(H, X,Y)$
on classical generators $(J_0,J_{+},J_{-})$. The full coalgebra structure 
of ${\cal U}_h(sl(2))$ can also be implemented. 

\smallskip

We have presented before [8,9] an alternative realization of the 
${\cal U}_h(sl(2))$ generators in terms of the classical generators:
\begin{eqnarray}
e^{h{\tilde X}}={I+{h\over 2}J_{+} \over I-{h\over 2}J_{+}} \qquad \Rightarrow
 \qquad 
{h{\tilde X} \over 2}=\hbox{arctanh}({hJ_{+} \over 2}),
 \end{eqnarray} 
whereas the remaining maps read
\begin{eqnarray}
&& {\tilde Y}=(1-(hJ_{+}/2)^2)^{1/2}J_{-}(1-(hJ_{+}/2)^2)^{1/2},
\qquad\qquad\qquad {\tilde H}=J_0.
 \end{eqnarray} 
Both the sets $({\hat H},{\hat X},{\hat Y})$ and $({\tilde H},{\tilde X},
{\tilde Y})$ satisfy the same ${\cal U}_h(sl(2))$ algebra. Strictly speaking 
such realizations hold on finite dimensional vectors spaces of irreducible 
representations. The nilpotency of the classical $J_+$ then, via the 
construction of the triplets $({\hat H},{\hat X},{\hat Y})$ and $({\tilde H},
{\tilde X},{\tilde Y})$, provides finite dimensional irreducible 
representations of the ${\cal U}_h(sl(2))$ algebra.  
 Inverting the maps (26) and (31) we obtain
\begin{eqnarray}
J_{+}= {1 \over h} \sinh(h{\hat X}) ={2 \over h} \tanh({h{\tilde X} \over 2}) .
 \end{eqnarray} 
In (29), we obtain that ${\hat Y}$ is nonlinear in $J_0$. This is
a consequence of the nonlinearity involved in the definition (27) of 
${\hat H}$. 

\smallskip

So far, as a relatively simple, illustrative example of our method, we 
have been considering the case $({1\over 2}\otimes j)$. But our method can 
be used to obtain the universal ${\cal R}_h$-matrix [4,5] in the form
   \begin{eqnarray}    
{\cal R}_h=\exp\biggl(-h{\hat X} \otimes e^{h{\hat X}}{\hat H}\biggr) 
\exp\biggl(he^{h{\hat X}}{\hat H}\otimes {\hat X}\biggr). 
 \end{eqnarray} 
In obtaining this form as the $q\rightarrow 1$ limit of transformed universal
${\cal R}$-matrix of ${\cal U}_q(sl(2))$ the cancellation of divergences
 become more subtle and complicated. We have fully analyzed it for the case 
$(1\otimes j)$. The procedure will be presented elsewhere. But (34) must 
evidently have all the required properties. 

\smallskip

We can study representations of ${\cal U}_q(sl(2))$ using the map
$({\hat H},{\hat X}, {\hat Y})\longrightarrow (J_0,J_{\pm})$
as it was done in [8,9] using $({\tilde H},{\tilde X},{\tilde Y})
\longrightarrow (J_0,J_{\pm})$.
This aspect will not be treated here. For the particularly simple example of 
$({1\over 2}\otimes 1)$ representation, we compare
\begin{eqnarray}    
R_{h}({\hat H},{\hat X},{\hat Y})=\pmatrix{{\hat A} && {\hat B} \cr
&& \cr
0 && {\hat C}\cr}
 \end{eqnarray}
with
\begin{eqnarray}    
R_{h}({\tilde H},{\tilde X},{\tilde Y}
)=\pmatrix{{\tilde A} && {\tilde B} \cr
&& \cr
0 && {\tilde C}\cr}.
 \end{eqnarray}
For the representation $({1\over 2}\otimes 1)$, 
$R_{h}({\tilde H},{\tilde X},{\tilde Y})$
was already given in [8]:
\begin{eqnarray}    
{\tilde A}=\pmatrix{1 &2h & 2h^2 \cr
0&1&2h \cr
0 &0& 1\cr}, \qquad {\tilde B}=\pmatrix{-2h &2h^2 & 0 \cr
0&0&2h^2 \cr
0 &0& 2h\cr},\qquad {\tilde C}=\pmatrix{1 &-2h & 2h^2 \cr
0&1&-2h \cr
0 &0& 1\cr}.
 \end{eqnarray}
Using our new mapping, we obtain 
\begin{eqnarray}    
{\hat A}={\tilde A},\qquad {\hat B}=\pmatrix{-2h &2h^2 & 4h^3 \cr
0&0&2h^2 \cr
0 &0& 2h\cr}, \qquad {\hat C}={\tilde C}.
\end{eqnarray}
The two $R$-matrices (35) and (36) are related through a similarity 
transformation by a matrix of the form 
$\pmatrix{1 & 0 \cr
0 & 1\cr} \otimes M$, where
\begin{eqnarray}    
M=\pmatrix{1 &0 & h^2 \cr
0&1&0 \cr
0 &0& 1\cr}.
\end{eqnarray}

\smallskip

The major interest in our method is that {\em it can be generalized for 
obtaining the nonstandard $R$-matrices of algebras of higher dimensions as 
contraction limits of the corresponding $R_q$-matrices}.  Here we illustrate 
the scope of our method, using $sl(3)$ as an example. We treat the case 
of $R$-matrices corresponding to the representation 
{\em (fund.)$\otimes$(arb.irrep.)}. Possibilities of further generalizations 
are then indicated. Here we assume the ${\cal U}_q(sl(3))$ Hopf algebra to be 
well-known [6]. Let the standard Chevally generators of ${\cal U}_q(sl(3))$ 
be $(q^{\pm h_1},q^{\pm h_2},{\hat e_1},{\hat e_2},{\hat f_1},{\hat f_2})$, 
whereas the classical $(q=1)$ generators are denoted by
$(h_1,h_2,e_1,e_2,$ $f_1,f_2)$. Defining
\begin{eqnarray}    
{\hat e}_3={\hat e}_1{\hat e}_2-q^{-1}{\hat e}_2{\hat e}_1, \qquad\qquad
 {\hat f}_3={\hat f}_2{\hat f}_1-q{\hat f}_1{\hat f}_2, 
\end{eqnarray}
we obtain
\begin{eqnarray}    
&& q^{h_i} {\hat e}_3={\hat e}_3q^{h_i+1}\;\;\hbox{for}\;\;(i=1,2),
\qquad\qquad[{\hat e}_3,{\hat f}_3]=[h_1+h_2],\nonumber \\
&&[{\hat e}_3,{\hat f}_1]=-{\hat e}_2q^{-h_1},\qquad\;\;\;\;\qquad\qquad\qquad
[{\hat e}_3,{\hat f}_2]=q^{h_2}{\hat e}_1,\nonumber \\
&& {\hat e}_1{\hat e}_3=q{\hat e}_3{\hat e}_1, 
\qquad\qquad\;\;\;\;\;\qquad \qquad\qquad
{\hat e}_2{\hat e}_3=q^{-1}{\hat e}_3{\hat e}_2 .
\end{eqnarray}
 As was noted in [3], after using the contraction process for ${\cal U}_q
(sl(2))$ subalgebra and possible relabeling of the generators, the really 
new possibility for ${\cal U}_q(sl(3))$ arises from using ${\hat e_3}$ in
the construction of the operator performing the similarity transformation. 
More generally, similar new possibility for the ${\cal U}_q(sl(N))$ 
emerges by using the generators corresponding to the highest positive root 
for the same purpose. As a direct generalization of (6) we now consider 
 \begin{eqnarray}
E_q(\eta {\hat e}_{3})=\sum_{n= 0}^{\infty} {(\eta {\hat e}_{3})^n\over [n]!}.
\end{eqnarray}
Using (41), we obtain  
\begin{eqnarray}    
&& {\hat e}_1E_q(\eta {\hat e}_3) =E_q(q\eta  {\hat e}_3){\hat e}_1, 
\qquad \qquad \qquad \;\;\qquad {\hat e}_2E_q(\eta {\hat e}_3) 
=E_q(q^{-1}\eta  {\hat e}_3){\hat e}_2,\nonumber \\
&& [E_q(\eta {\hat e}_3) ,f_1]=-\eta E_q(q^{-1}\eta {\hat e}_3)
{\hat e}_2q^{-h_1}\qquad\qquad 
[E_q(\eta{\hat e}_3),f_2]=\eta E_q(q\eta{\hat e}_3)q^{h_2}{\hat e}_1\nonumber\\
&& [E_q(\eta {\hat e}_3),f_3]={\eta \over q-q^{-1}}
\biggl(E_q(q\eta {\hat e}_3)q^{h_1+h_2}-
E_q(q^{-1}\eta {\hat e}_3) q^{-h_1-h_2}\biggr)
\end{eqnarray}
Analogous to the example of ${\cal U}_q(sl(2))$ algebra, we here define 
\begin{eqnarray}
&& {\hat {\tt T}}_{(\alpha)}=(E_q(\eta {\hat e_{3}}))^{-1}E_q(q^{\alpha}\eta 
{\hat e_{3}}), \\
&& {\tt T}_{(\alpha)}=\lim_{q\rightarrow 1} {\hat {\tt T}}_{(\alpha)}.
\end{eqnarray}
To obtain the finite limit of the operators ${\tt T}_{(\alpha)}$ we may 
proceed by exploiting the ${\cal U}_q(sl(2))$ subalgebra generated by 
$(e_3,f_3,q^{\pm(h_1+h_2)})$. Equivalently, we note that the {\em in the 
formal series development of $E_q(x)$, as 
well as in the subsequent limiting process no specific properties of 
${\cal J}_{+}$ (or ${\hat e_3}$) need to be introduced}. Parallel to (22) we 
here obtain   
\begin{eqnarray}
&&{\tt T}_{(\pm 1)}=\pm he_3+\biggl(1+(he_3)^2\biggr)^{1/2},\qquad\qquad
\qquad{\tt T}_{(\alpha)}={\tt T}^{\alpha}.
\end{eqnarray}
For ${\cal U}_q(sl(3))$ algebra the $R_q$ matrix in the representation 
{\em (fund.)$\otimes$(arbitrary irrep.)} reads [6]:
\begin{eqnarray}
R_q=\pmatrix{q^{{1\over 3}(2h_1+h_2)} &q^{{1\over 3}(2h_1+h_2)}\Lambda_{12}
 &q^{{1\over 3}(2h_1+h_2)}\Lambda_{13} \cr 
0 &q^{-{1\over 3}(h_1-h_2)}
 &q^{-{1\over 3}(h_1-h_2)}\Lambda_{23} \cr
 0 & 0
 &q^{-{1\over 3}(h_1+2h_2)}\cr}, 
\end{eqnarray}
where
\begin{eqnarray}
&& \Lambda_{12}=q^{-1/2}(q-q^{-1})q^{-h_1/2}{\hat f_1},\qquad\qquad
 \Lambda_{23}=q^{-1/2}(q-q^{-1})q^{-h_2/2}{\hat f_2}, \nonumber\\
&& \qquad\qquad\qquad\qquad \Lambda_{13}=q^{-1/2}(q-q^{-1}){\hat f_3}
q^{-{1\over 2}(h_1+h_2)}.
\end{eqnarray}
We introduce the operator  
\begin{eqnarray}
{\tt G}=E_q(\eta {\hat e_3})_{(fund.)} \otimes E_q(\eta {\hat e_3})_{(arb.)} 
\end{eqnarray}
and perform similarity transformation on the $R_q$ matrix (47) 
 \begin{eqnarray}
&& {\tt G}^{-1}R_q{\tt G}= 
\pmatrix{{\tt g}^{-1}q^{{1\over 3}(2h_1+h_2)}{\tt g} 
&{\tt g}^{-1}q^{{1\over 3}(2h_1+h_2)}\Lambda_{12}{\tt g}
 &\eta\;{\tt g}^{-1}\bigl(q^{{1\over 3}(2h_1+h_2)}-q^{-{1\over 3}(h_1+2h_2)}
\bigr){\tt g} \cr
& & +{\tt g}^{-1}q^{{1\over 3}(2h_1+h_2)}\Lambda_{13}{\tt g}\cr
&& \cr 
0 &{\tt g}^{-1}q^{-{1\over 3}(h_1-h_2)}{\tt g}
 &{\tt g}^{-1}q^{-{1\over 3}(h_1-h_2)}\Lambda_{23}{\tt g} \cr
&& \cr && \cr
 0 & 0 &{\tt g}^{-1}q^{-{1\over 3}(h_1+2 h_2)}{\tt g}\cr},\;\;\; 
\end{eqnarray}
where we use the notation ${\tt g}=E_q(\eta {\hat e_3})$. In close analogy 
with earlier example, the transform turns out to be finite at the 
$q\rightarrow 1$ limit: 
\begin{eqnarray}
\lim_{q\rightarrow 1} {\tt G}^{-1}R_q{\tt G}=R_h=\pmatrix{T & 2hT^{-1/2}e_2
 & -{h\over 2}(T+T^{-1})(h_1+h_2)+{h \over 2}(T-T^{-1}) \cr 
0 & I & -2hT^{1/2}e_1\cr
 0 & 0
 & T^{-1}\cr}. 
\end{eqnarray}
We have, briefly, directly stated the final result. But the derivation is 
fairly analogous to $sl(2)$ case. {\em In fact this illustrates how directly 
our formalism can be implemented for higher dimensional cases}. 

\smallskip

We note that the corner elements of (51) have exactly the same 
structure as (24). The new features arises with the presence of terms 
$2hT^{-1/2} e_2$, $-2hT^{1/2} e_1$ involving the simple root generators. 
In general, we may proceed as follows. The matrix structure in (51) provides 
a realization of the Lax operator corresponding to the $h$-Borel subalgebra 
${\cal U}_h({\cal B}_+)$ generated by the Cartan elements and the positive
root generators. In fact following the analogy for the ${\cal U}_h(sl(2))$
case, as evidenced by comparing (1) with (25), the matrix operator (51) 
provides a map of the $h$-Borel subalgebra ${\cal U}_h({\cal B}_+)$ on its 
classical counterpart. The full Hopf structure of the $h$-Borel subalgebra
may then be obtained by the standard FRT procedure [10]. In the instance of 
nonstandard $h$-deformed $sl(2)$ algebra, the universal $R_h$-matrix for 
this Borel subalgebra is the universal $R_h$-matrix for the full 
${\cal U}_h(sl(2))$ Hopf algebra. Assuming that this property still holds 
for all $h$-deformed $sl(N)$ algebra, then our construction may provide a 
route to obtain the so far uninvestigated ${\cal U}_h(sl(N))$ algebras 
for $N>2$. In this instance, it is worth pointing out the close kinship
of the $h$-deformed algebras with $\kappa$-Poincar\'e and deformed conformal 
algebras [11], where the deformation parameter has dimensions of mass and,
consequently, an induced fundamental length scale enters at the geometrical 
level.
         
\smallskip

To elucidate a further possible application of our present technique of 
extraction of $h$-deformed objects from the corresponding $q$-deformed 
objects via contraction procedure, we return to the earlier example of 
deformed $sl(2)$ algebra. Exploiting the bialgebraic duality relationship
between the pair of Hopf algebras ${\cal F}un_q(GL(2))$ and 
${\cal U}_q(gl(2))$, 
the dual form, which provides a generalization of the familiar exponential 
relationship between classical Lie groups and algebras, was constructed in 
[12]. This may be used to provide [13] arbitrary finite dimensional irreducible
representations of $GL_q(2)$ leading to the $q$-analogues of the classical 
Wigner $d^{(j)}$-functions or spherical functions. The present method 
of using similarity transformation, singular itself at $q\rightarrow 1$ limit, 
but in such a way that the singularities get cancelled to 
provide finite results for the transformed objects, may then be used to obtain 
arbitrary finite dimensional representations of the $h$-deformed group 
elements $GL_h(2)$ and, consequently, $h$-deformed analogues of Wigner 
$d^{(j)}$ functions. We will return to this topic later.

\smallskip

At the end we wish to mention a fact as a curiosity. There exist canonical 
transformations of the classical Heisenberg algebra closely mimicking our 
construction of the nonlinear maps of ${\cal U}_h(sl(2))$ generators on
the generators of the classical $sl(2)$ algebra. Starting with classical 
Heisenberg commutation relation $[a_-,a_+]=1$, we define an one-parameter 
$(h)$ family of transformations 
\begin{eqnarray}      
{\hat a_{+}}={1\over h}\;\ln\bigl(ha_{+}+\sqrt{1+h^2a_{+}^2}\bigr),
\qquad\qquad\qquad {\hat a_{-}} =  \sqrt{1+h^2a_{+}^2}\;\;a_{-},
\end{eqnarray}
which satisfy Heisenberg commutation relation $[{\hat a_-},{\hat a_+}]=1$.
This closely resembles the nonlinear realization of the generators of 
${\cal U}_h(sl(2))$ discussed in (26), (27) and (28). There exists also an 
anlogous family of canonical transformations closely paralleling the other 
nonlinear realization of ${\cal U}_h(sl(2))$ generators discussed in (31) 
and (32):
\begin{eqnarray}      
{\tilde a_{+}}={2\over h}\;\hbox{arctanh}({ha_{+}\over 2}),
\qquad\qquad\qquad {\tilde a_{-}} = \sqrt{1-{h^2a_{+}^2 \over 4}}
\;\;a_{-}\;\;\sqrt{1-{h^2a_{+}^2 \over 4}},
\end{eqnarray}
which again satisfies the commutation relation $[{\tilde a_-},{\tilde a_+}]=1$.
It is perhaps worth investigating the possible connection between the theory
of canonical transformations of the Heisenberg algebra and the nonlinear map
discussed above.

\vskip 1cm

\noindent {\bf Acknowledgments:} 
 
We thank R. Jagannathan for interesting discussions. 
One of us (RC) wants to thank A. Chakrabarti for a kind invitation. He is 
also grateful to the members of the CPTH group for their kind hospitality.

\end{document}